\documentclass[11pt,reqno,sumlimits]{amsart}
\usepackage{amsfonts, amsmath, amscd, amssymb, euscript, amsthm, array, booktabs,
dcolumn, shortvrb, tabularx, units, url, mathrsfs}
\usepackage[pdftex]{graphicx}
\usepackage[pdftex]{hyperref}
\usepackage[all]{xy}
\normalsize
\DeclareMathOperator{\cl}{C\ell}

\DeclareMathOperator{\Aut}{Aut}
\DeclareMathOperator{\fl}{flat}
\DeclareMathOperator{\LL}{L}
\DeclareMathOperator{\Cl}{C\ell}

\DeclareMathOperator{\topo}{t}
\DeclareMathOperator{\todd}{Todd}
\DeclareMathOperator{\an}{a}

\DeclareMathOperator{\ind}{ind}

\DeclareMathOperator{\spin}{spin}

\DeclareMathOperator{\CS}{CS}

\DeclareMathOperator{\odd}{odd}
\DeclareMathOperator{\even}{even}

\DeclareMathOperator{\im}{Im}
\DeclareMathOperator{\id}{id}
\DeclareMathOperator{\ch}{ch}

\DeclareMathOperator{\End}{End}
\DeclareMathOperator{\ho}{Hom}
\DeclareMathOperator{\rk}{rank}
\DeclareMathOperator{\str}{str}
\DeclareMathOperator{\tr}{tr}
\begin{document}
\setlength{\baselineskip}{1.5\baselineskip}
\theoremstyle{definition}
\newtheorem{coro}{Corollary}
\newtheorem*{note}{Note}
\newtheorem{thm}{Theorem}
\newtheorem*{ax}{Axiom}
\newtheorem{defi}{Definition}
\newtheorem{lemma}{Lemma}
\newtheorem*{claim}{Claim}
\newtheorem{exam}{Example}
\newtheorem{prop}{Proposition}
\newtheorem{remark}{Remark}
\newcommand{\wt}[1]{{\widetilde{#1}}}
\newcommand{\ov}[1]{{\overline{#1}}}
\newcommand{\wh}[1]{{\widehat{#1}}}
\newcommand{\poin}{Poincar\'e~}
\newcommand{\deff}[1]{{\bf\emph{#1}}}
\newcommand{\boo}[1]{\boldsymbol{#1}}
\newcommand{\abs}[1]{\lvert#1\rvert}
\newcommand{\norm}[1]{\lVert#1\rVert}
\newcommand{\inner}[1]{\langle#1\rangle}
\newcommand{\poisson}[1]{\{#1\}}
\newcommand{\biginner}[1]{\Big\langle#1\Big\rangle}
\newcommand{\set}[1]{\{#1\}}
\newcommand{\Bigset}[1]{\Big\{#1\Big\}}
\newcommand{\BBigset}[1]{\bigg\{#1\bigg\}}
\newcommand{\dis}[1]{$\displaystyle#1$}
\newcommand{\R}{\mathbb{R}}
\newcommand{\EE}{\mathbb{E}}
\newcommand{\GG}{\mathbb{G}}
\newcommand{\N}{\mathbb{N}}
\newcommand{\Z}{\mathbb{Z}}
\newcommand{\Q}{\mathbb{Q}}
\newcommand{\E}{\mathcal{E}}
\newcommand{\G}{\mathcal{G}}
\newcommand{\F}{\mathcal{F}}
\newcommand{\V}{\mathcal{V}}
\newcommand{\W}{\mathcal{W}}
\newcommand{\SSS}{\mathcal{S}}
\newcommand{\h}{\mathbb{H}}
\newcommand{\C}{\mathbb{C}}
\newcommand{\A}{\mathcal{A}}
\newcommand{\MM}{\mathcal{M}}
\newcommand{\HH}{\mathcal{H}}
\newcommand{\D}{\mathcal{D}}
\newcommand{\PP}{\mathcal{P}}
\newcommand{\M}{\mathcal{M}}
\newcommand{\K}{\mathcal{K}}
\newcommand{\RRR}{\mathscr{R}}
\newcommand{\AAA}{\mathscr{A}}
\newcommand{\DDD}{\mathscr{D}}
\newcommand{\e}{\mathscr{E}}
\newcommand{\f}{\mathscr{F}}
\newcommand{\g}{\mathscr{G}}
\newcommand{\so}{\mathfrak{so}}
\newcommand{\gl}{\mathfrak{gl}}
\newcommand{\aaa}{\mathbb{A}}
\newcommand{\bbb}{\mathbb{B}}
\newcommand{\DD}{\mathsf{D}}
\newcommand{\FF}{\mathbb{F}}
\newcommand{\ccc}{\bold{c}}
\newcommand{\sss}{\mathbb{S}}
\newcommand{\cdd}[1]{\[\begin{CD}#1\end{CD}\]}
\numberwithin{equation}{subsection}
\normalsize
\title[Flat GRR]{The flat Grothendieck--Riemann--Roch theorem without adiabatic
techniques}
\author[M.-H. Ho]{Man-Ho Ho}
\address{Department of Mathematics\\ Hong Kong Baptist University}
\email{homanho@math.hkbu.edu.hk}
\maketitle
\nocite{*}
\begin{abstract}
In this paper we give a simplified proof of the flat Grothendieck--Riemann--Roch
theorem. The proof makes use of the local family index theorem and basic computations
of the Chern--Simons form. In particular, it does not involve any adiabatic limit
computation of the reduced eta-invariant.
\end{abstract}
\tableofcontents
\section{Introduction}

In this paper we give a simplified proof of the flat Grothendieck--Riemann--Roch theorem
that avoids adiabatic techniques.

\subsection{Historical background}

In this subsection we briefly review the historical background of flat $K$-theory, its
Grothendieck--Riemann--Roch theorem and its relation to physics. For a detailed
exposition, see \cite{BS12, F00}.

$\R/\Z$ $K$-theory \cite{APS76} (also called flat $K$-theory) is a generalized cohomology
theory which is defined as the cokernel of a natural homomorphism $K^*(X; \Q)\to K^*(X;
\Q/\Z)\oplus K^*(X; \R)$. One motivation of defining $\R/\Z$ $K$-theory is to prove a
cohomological version of the Atiyah--Singer family index theorem \cite{AS71} (FIT) for
bundles with vanishing Chern characters: if $[E]\in K(X)$ has vanishing Chern character,
the Grothendieck--Riemann--Roch theorem (GRR)
\begin{equation}\label{eq 1.1.1}
\ch(\ind^{\an}(E))=\int_{X/B}\todd(X/B)\cup\ch(E)
\end{equation}
implies that $\ch(\ind^{\an}(E))=0$. (Here $X\to B$ is a fibration of closed manifolds,
$E\to X$ is a complex vector bundle, and all other terms are defined in the paper.) Thus
one should get a refinement of GRR for bundles with vanishing Chern characters. It turns
out that the flat $K$-group $K^{-1}(X; \R/\Z)$ is the right home for such a refinement.

The first geometric model of $K^{-1}(X; \C/\Z)$ is given by Karoubi \cite{Kar87} under
the name ``multiplicative $K$-theory". By adding Hermitian structures to elements in
$K^{-1}(X; \C/\Z)$, Lott gives the geometric model of $K^{-1}(X; \R/\Z)$, which is
denoted by $K^{-1}_{\LL}(X)$ in this paper, and proves the FIT in $K^{-1}_{\LL}$
\cite{L94} (flat FIT), which equates the flat analytic index $\ind^{\an}_{\LL}$ and the
topological index $\ind^{\topo}$. The flat FIT is the refinement of the FIT for bundles
with vanishing Chern characters. The corresponding GRR \cite[Corollary 4]{L94} (flat GRR)
\begin{equation}\label{eq 1.1.2}
\ch_{\R/\Q}(\ind^{\an}_{\LL}(\E))=\int_{X/B}\todd(X/B)\cup\ch_{\R/\Q}(\E)\in H^{\odd}
(B; \R/\Q),
\end{equation}
where $\ch_{\R/\Q}:K^{-1}_{\LL}(B)\to H^{\odd}(B; \R/\Q)$ is the flat Chern character
and $\ind^{\an}_{\LL}:K^{-1}_{\LL}(X)\to K^{-1}_{\LL}(B)$ is the flat analytic index,
is the refinement of (\ref{eq 1.1.1}).

In modern language $K^{-1}_{\LL}$ is the flat part of differential $K$-theory $\wh{K}^0$
constructed by Hopkins--Singer \cite{HS05}, Bunke--Schick \cite{BS09}, Freed--Lott
\cite{FL10} and Simons--Sullivan \cite{SS10} respectively, which is a generalized
differential cohomology theory in the sense of Bunke--Schick \cite{BS10}. In theoretical
physics the motivation of differential $K$-theory comes from the assertion of Witten
\cite{W98} that D-brane charges in string theory are described by a $K$-theory class of
spacetime rather than by a cohomology class. Furthermore, Moore and Witten propose that
Ramond--Ramond fields in type II and type I string theory, to which D-branes couple, are
also classified by $K$-theory \cite{MW00}. Freed and Hopkins propose using differential
$K$-theory to describe Ramond--Ramond fields \cite{DH00}.

Freed and Lott prove a FIT in differential $K$-theory \cite[Theorem 7.35]{FL10} (dFIT),
which equates the differential analytic index and the differential topological index.
The GRR in differential $K$-theory \cite[Corollary 8.26]{FL10} (dGRR) is also
proved by Bunke and Schick \cite[Theorem 6.19]{BS09} independently. The flat FIT and
the flat GRR can be considered as special cases of the dFIT and dGRR respectively. See
\cite{HS15} for an algebraic analog of differential cohomology and the corresponding
Riemann--Roch theorem. The motivation in theoretical physics for formulating and proving
the dFIT, or rather its consequence for determinant line bundles dates back for proving
the Green--Schwarz cancellation of local and global anomalies in type I string theory
\cite{F00}.

\subsection{Outline of proof and the relation to previous work}

In this subsection we first outline our proof of the flat GRR. Then we discuss the
relation between our proof and the previous proofs and raise some questions.

First of all we briefly outline our proof of the flat GRR. Let $\pi:X\to B$ be a
submersion with closed $\spin^c$ fibers of even relative dimension. Consider the
associated submersion $\pi\times\id:X\times I\to B\times I$, where $I=[0, 1]$. The
local FIT of the $\spin^c$ Dirac operator twisted by a $\Z_2$-graded Hermitian bundle
$\e\to X\times I$ with a $\Z_2$-graded unitary connection $\nabla^{\e}$ is given by
\begin{equation}\label{eq 1.2.1}
d\wt{\eta}(\e)=\int_{X\times I/B\times I}\todd(\nabla^{S^c(T^V(X\times I))})\wedge\ch
(\nabla^{\e})-\ch(\nabla^{\ker(\DD^{\e})}),
\end{equation}
where $\wt{\eta}(\e)$ is the Bismut--Cheeger eta form \cite{BC89, D91}. Here we have
assumed that the family of the complex vector spaces $\ker(\DD^{\e}_z)$ form a vector
bundle $\ker(\DD^{\e})\to B\times I$. Note that (\ref{eq 1.2.1}) is an equality of
differential forms on $B\times I$. By integrating (\ref{eq 1.2.1}) along the fibers of
the trivial fibration $B\times I\to B$ we obtain the variational formula of the eta
forms. We also prove the additivity of the eta forms. These two results enable us to
prove that the flat analytic index is well defined (Proposition \ref{prop 3}). To prove
the flat GRR we choose a suitable $\e\to X\times I$, which is the pullback of a certain
$\Z_2$-graded Hermitian bundle over $X$, in the variational formula of the eta forms.
This will give us an equality of closed differential forms of odd degree on $B$, whose
mod $\Q$ reduction of its de Rham class is (\ref{eq 1.1.2}) (Theorem \ref{thm 1}). For
the general case where the family of the complex vector spaces $\ker(\DD^{\e}_z)$ do
not form a vector bundle, one can prove the corresponding results along the lines of
\cite[\S5]{L94} (see also \cite[\S7]{FL10}). All our arguments in the special case carry
over to the general case.

One important ingredient of the previous proofs is the adiabatic limit of the reduced
eta-invariant of spin (or $\spin^c$) Dirac operator, which we briefly recall. Let $B$
be a closed odd-dimensional spin manifold and $\pi:X\to B$ a submersion with closed
spin fibers of even relative dimension. For $\varepsilon>0$, consider the submersion
metric $g^{TX}_\varepsilon=\varepsilon^{-1}\pi^*g^{TB}\oplus g^{T^VX}$ with respect to
a horizontal distribution. Denote by $\DD_\varepsilon$ the corresponding spin Dirac
operator and by $\bar{\eta}(\DD_\varepsilon)$ the associated reduced eta-invariant. The
study of the limiting behaviour, for example, of $\bar{\eta}(\DD_\varepsilon)$ when
$\varepsilon\to 0$, is called passing to adiabatic limit. It is rooted in \cite{AS84}
and initiated by Witten \cite{W85}, who relates the adiabatic limit of the eta-invariant
to the holonomy of determinant line bundle, the global anomaly. Witten's result receives
rigourous proofs in \cite{BF86a, BF86b, C87, DF94}. Adiabatic limit becomes an important
tool in, among many other areas, local index theory (see \cite{BC89, D91, BGV} and the
references therein).

Now we briefly outline the previous proofs of the flat GRR. As mentioned in \S1.1 the
flat GRR is a direct consequence of the flat FIT. One could also prove the flat GRR
directly in the spirit of \cite{L94}, which shares some similarities to the proof of
the flat FIT and the proof of \cite[Proposition 8.19]{FL10}.

The universal coefficient theorem for ordinary cohomology and the divisibility of
$\R/\Q$ imply that $H^{\odd}(B; \R/\Q)\cong\ho(H_{\odd}(B); \R/\Q)$, or equivalently the
existence of a pairing $\inner{, }_H:H^{\odd}(B; \R/\Q)\times H_{\odd}(B)\to\R/\Q$. As
(\ref{eq 1.1.2}) is an equality in $H^{\odd}(B; \R/\Q)$, proving it is equivalent to
proving
\begin{equation}\label{eq 1.2.2}
\biginner{\ch_{\R/\Q}(\ind^{\an}_{\LL}(\E))-\int_{X/B}\todd(X/B)\cup\ch_{\R/\Q}
(\E), U}_H
\end{equation}
is zero in $\R/\Q$ for every $U\in H_{\odd}(B)$. There is a pairing analogous to
$\inner{, }_H$ on the $K$-theory level, guaranteed by the universal coefficient
theorem for generalized cohomology theory \cite{Y75} and the divisibility of $\R/\Z$.
Denote by $\inner{, }_K:K^{-1}_{\LL}(B)\times K_{-1}(B)\to\R/\Z$ the pairing, where
$K_\bullet$ is the topological $K$-homology group given by Baum--Douglas \cite{BD82}.
The pairings $\inner{, }_H$ and $\inner{, }_K$ are related by the flat Chern character
$\ch_{\R/\Q}$ and the homological Chern character $\ch_{\odd}:K_{-1}(B)\to H_{\odd}
(B; \Q)$ in the sense that the following diagram commutes.
\begin{equation}\label{eq 1.2.3}
\begin{CD}
K^{-1}_{\LL}(B)\times K_{-1}(B) @>\inner{, }_K>> \R/\Z \\ @V\ch_{\R/\Q}\times
\ch_{\odd}VV @VVV \\ H^{\odd}(B; \R/\Q)\times H_{\odd}(B; \Q) @>>\inner{, }_H> \R/\Q
\end{CD}
\end{equation}
The pairing $\inner{, }_K$ can be given by the reduced eta-invariant
\cite[Proposition 3]{L94} as follows. For a $\Z_2$-graded generator $\E=(E^+\oplus E^-,
h^+\oplus h^-, \nabla^+\oplus\nabla^-, \omega)$ of $K^{-1}_{\LL}(B)$ (see \S2.2 for the
details) and a cycle $\K=(X, F, f)$ of $K_{-1}(B)$, which consists of a complex vector
bundle $F\to X$ over a closed odd-dimensional $\spin^c$ manifold and a smooth map $f:X
\to B$, define $\bar{\eta}(f^*\E)\in\R/\Z$ by
\begin{equation}\label{eq 1.2.4}
\bar{\eta}(f^*\E):=\bar{\eta}(\DD^{F\otimes f^*E^+})-\bar{\eta}(\DD^{F\otimes f^*E^-})-
\int_X\todd(\nabla^{S^c(TX)})\wedge\ch(\nabla^F)\wedge f^*\omega.
\end{equation}
Then
\begin{equation}\label{eq 1.2.5}
\inner{[\E], [\K]}_K=\bar{\eta}(f^*\E).
\end{equation}
Since $\ch_\bullet:K_\bullet(B)\otimes\Q\to H_\bullet(B; \Q)$ is an isomorphism (in
particular it is surjective), it follows from the arguments in the proof of
\cite[Proposition 8.19]{FL10} and \cite[Proposition 6]{L94} that one can take $U$ in
(\ref{eq 1.2.2}) to be $\ch_{\odd}([B])$, where $B$ is now assumed to be a closed
odd-dimensional $\spin^c$-manifold and $[B]$ is the fundamental $K$-homology class. By
(\ref{eq 1.2.3}) and (\ref{eq 1.2.5}), proving (\ref{eq 1.2.2}) is zero in $\R/\Q$ boils
down to computing the reduced eta-invariants of some $\spin^c$ Dirac operators and its
adiabatic limits.

On the other hand, one can apply \cite[Theorem 1.15]{B05} to prove the flat GRR, which
is done in the previous version of this paper. The proof of \cite[Theorem 1.15]{B05} is
somewhat similar to the above proof, as it also consists of computations of the reduced
eta-invariants of spin Dirac operators and their adiabatic limits.

One might ask if there is any relation between our proof of the flat GRR and the previous
proofs. Since our proof of the flat GRR does not involve any adiabatic limit of the
reduced eta-invariant, one might wonder whether the some of the results in
\cite{BF86a, BF86b, C87, BC89, D91} can be proved without it. Frankly we do not have any
informative answers for these questions at this moment. Perhaps a clue for these
questions can be found in \cite[\S4]{APS75b}, which is an interesting topic to be further
investigated.

Since the flat GRR is a special case of the dGRR, one suspects that whether
\cite[Theorem 1.15]{B05} or even the dFIT can be proved without computing adiabatic limit
of the reduced eta-invariant. Note that \cite[Theorem 1.15]{B05} takes values in
Cheeger--Simons differential characters \cite{CS85}. Our experience shows that equality
of differential characters is usually harder to prove than equality of differential forms.
More precisely, the proofs of \cite[Theorem 1.15]{B05}, the dGRR and the dFIT depend
crucially on \cite[Theorem 9.2]{CS85} and \cite[Theorem 0.1']{D91} (see also
\cite[(52)]{L94}). Thus the affirmative answer to this question depends on the previous
questions.

This paper is organized as follows. In Section 2 we review the background material,
including some aspects of Chern--Weil theory, the flat $K$-theory, the setup and the
statement of the local FIT, and the definition of the flat analytic index. In Section 3
we prove the main results in this paper.

\section*{Acknowledgement}

The author would like to thank the referee for raising some interesting questions and
providing invaluable comments and suggestions which improve the paper in a significant
way.

\section{Background material}

In this paper $X$ and $B$ are closed manifolds and $I$ is the closed interval $[0, 1]$.
Given a manifold $X$, write $\wt{X}=X\times I$. Given $t\in[0, 1]$, define two maps
$i_{X, t}:X\to\wt{X}$ and $p_X:\wt{X}\to X$ by $i_{X, t}(x)=(x, t)$ and $p_X(x, t)=x$.

\subsection{Chern character form and Chern--Simons form}

Let $E\to X$ be a complex vector bundle with a Hermitian metric $h^E$ and a unitary
connection $\nabla^E$. The Chern character form of $\nabla^E$ is defined by
$$\ch(\nabla^E)=\tr(e^{-\frac{1}{2\pi i}(\nabla^E)^2})\in\Omega^{\even}_\Q(X),$$
where $\Omega^{\even}_\Q(X)$ is the set of all closed even forms on $X$ with periods
in $\Q$.

There is a ``canonical" transgression form \dis{\CS(\nabla^E_1, \nabla^E_0)\in
\frac{\Omega^{\odd}(X)}{\im(d)}} between the Chern character forms of two
connections in the sense that
\begin{equation}\label{eq 2.1.1}
d\CS(\nabla^E_1, \nabla^E_0)=\ch(\nabla^E_1)-\ch(\nabla^E_0).
\end{equation}
Define $\CS(\nabla^E_1, \nabla^E_0)$ as follows. In the following $k\in\set{0, 1}$ is
fixed. Note that $p_X\circ i_{X, k}=\id_X$ and $i_{X, k}\circ p_X\sim\id_{\wt{X}}$. Let
$\e\to\wt{X}$ be a complex vector bundle with a Hermitian metric $h^{\e}$ and a unitary
connection $\nabla^{\e}$. Note that $\e\cong p_X^*(i_{X, k}^*\e)$. Thus
$$E_0:=i_{X, 0}^*\e\cong i_{X, 0}^*p_X^*(i_{X, 0}^*\e)\cong i_{X, 1}^*p_X^*(i_{X, 0}^*
\e)\cong i_{X, 1}^*\e=:E_1.$$
Write $E=E_0\cong E_1$. Define $h^{E_k}=i_k^*h^{\e}$. By
\cite[Corollary 8.9, Chapter 1]{K08} there exists $f\in\Aut(E)$ such that $h^{E_0}=f^*
h^{E_1}$, so we may assume that $h^{E_0}=h^{E_1}$ and denote it by $h^E$. Define
$$\nabla_k^E:=i_k^*\nabla^{\e}.$$
Note that the connection $\nabla^E_k$ is compatible with $h^{E_k}$. The assumption
$h^{E_0}=h^{E_1}$ implies that both $\nabla^E_0$ and $\nabla^E_1$ are compatible with
$h^E$. Define
\begin{equation}\label{eq 2.1.2}
\CS(\nabla^E_1, \nabla^E_0)=\int_{\wt{X}/X}\ch(\nabla^{\e})\mod\im(d),
\end{equation}
where $\wt{X}/X$ denotes the fiber of the fiber bundle $\wt{X}\to X$, and
\dis{\int_{\wt{X}/X}} denotes integration along the fiber.

To prove the Chern--Simons form defined by (\ref{eq 2.1.2}) satisfies (\ref{eq 2.1.1}),
we need to invoke Stokes' theorem for integration along the fibers
\cite[Problem 4 (p.331)]{GHV}. In general, for a smooth fiber bundle $M\to B$, where $M$
is a manifold with boundary, with compact fibers of dimension $n$ satisfying certain
orientability assumptions, we have
\begin{equation}\label{eq 2.1.3}
(-1)^{k-n}\int_{\partial M/B}i^*\omega=\int_{M/B}d_M\omega-d_B\int_{M/B}\omega,
\end{equation}
where $i:\partial M\to M$ is the inclusion map and $\omega\in\Omega^k(M)$. Applying
(\ref{eq 2.1.3}) to the fiber bundle $\wt{X}\to X$, we have
\begin{displaymath}
\begin{split}
d\CS(\nabla^E_1, \nabla^E_0)&=d\int_{\wt{X}/X}\ch(\nabla^{\e})=\int_{\wt{X}/X}d\ch
(\nabla^{\e})+\int_{\partial\wt{X}/X}i^*\ch(\nabla^{\e})\\
&=\ch(\nabla_1^E)-\ch(\nabla_0^E).
\end{split}
\end{displaymath}

Given a Hermitian bundle $E\to X$ with two unitary connections $\nabla^E_0$ and
$\nabla^E_1$, one can apply the above construction to $(p_X^*E, p_X^*h^E,
\nabla^{p_X^*E})$ with
\begin{equation}\label{eq 2.1.4}
\nabla^{p_X^*E}:=\nabla^E_t+dt\wedge\frac{\partial}{\partial t},
\end{equation}
where $\nabla^E_t$ is a smooth curve of unitary connections joining $\nabla^E_0$ and
$\nabla^E_1$. Note that $\CS(\nabla^E_1, \nabla^E_0)$ is independent of the choice of
$\nabla^E_t$ \cite[Proposition 1.1]{SS10}.

Another equivalent definition of the Chern--Simons form is given by
\begin{equation}\label{eq 2.1.5}
\CS(\nabla^E_1, \nabla^E_0)=\int^1_0\tr\bigg(\frac{d\nabla^E_t}{dt}e^{-\frac{1}{2\pi i}
(\nabla^E_t)^2}\bigg)dt\mod\im(d).
\end{equation}
It follows from (\ref{eq 2.1.2}) that the Chern--Simons form satisfies the following
properties:
\begin{eqnarray}
\CS(\nabla^E_1, \nabla^E_0)&=&-\CS(\nabla^E_0, \nabla^E_1)\label{eq 2.1.6},\\
\CS(\nabla^E_1, \nabla^E_0)&=&\CS(\nabla^E_1, \nabla^E_2)+\CS(\nabla^E_2,
\nabla^E_0),\label{eq 2.1.7}\\
\CS(\nabla^E_1\oplus\nabla^F_1, \nabla^E_0\oplus\nabla^F_0)&=&\CS(\nabla^E_1, \nabla^E_0)
+\CS(\nabla^F_1, \nabla^F_0),\label{eq 2.1.8}
\end{eqnarray}
where $\nabla^F_1, \nabla^F_0$ are unitary connections on the Hermitian bundle $F\to X$.
The proofs of (\ref{eq 2.1.6})-(\ref{eq 2.1.8}) using (\ref{eq 2.1.5}) are given in
\cite[Proposition 1.1, Lemma 1.4]{SS10}.

One can define the Chern character form and the Chern--Simons form of unitary
superconnection on $\Z_2$-graded Hermitian bundles in the exact same way as above, except
that the traces in the definitions are replaced by supertraces \cite{Q85},
\cite[\S1.4, \S1.5]{BGV}. Note that (\ref{eq 2.1.1}) and (\ref{eq 2.1.6})-(\ref{eq 2.1.8})
hold for unitary superconnections.

\subsection{The flat $K$-theory}

In this subsection we recall the flat $K$-theory \cite{L94}.

The flat $K$-group $K^{-1}_{\LL}(X)$ is an abelian group given by generators and relations:
a generator is of the form $\E=(E, h^E, \nabla^E, \omega)$, where \dis{\omega\in
\frac{\Omega^{\odd}(X)}{\im(d)}} satisfies $\ch(\nabla^E)-\rk(E)=-d\omega$.\footnote{This
differs from \cite[Definition 5]{L94} by a sign.} The only relation is $\E_1=\E_0$ if and
only if there exists $\G=(G, h^G, \nabla^G, \omega_G)$ such that $E_1\oplus G\cong E_0
\oplus G$ and
$$\omega_1-\omega_0=\CS(\nabla^E_0\oplus\nabla^G, \nabla^E_1\oplus\nabla^G).$$
Elements in $K^{-1}_{\LL}(X)$ are required to have virtual rank zero.

A $\Z_2$-graded generator $\E$ of $K^{-1}_{\LL}(X)$ has the form
\begin{equation}\label{eq 2.2.0}
\E=(E^+\oplus E^-, h^+\oplus h^-, \nabla^+\oplus\nabla^-, \omega),
\end{equation}
where $E^+\oplus E^-\to X$ is a $\Z_2$-graded complex vector bundle with a $\Z_2$-graded
Hermitian metric $h^+\oplus h^-$, a $\Z_2$-graded unitary connection $\nabla^+\oplus
\nabla^-$ on $E^+\oplus E^-\to X$, and \dis{\omega\in\frac{\Omega^{\odd}(X)}{\im(d)}}
satisfying
$$\ch(\nabla^+\oplus\nabla^-)=\ch(\nabla^+)-\ch(\nabla^-)=-d\omega.$$
Every element in $K^{-1}_{\LL}(X)$ can be written as a $\Z_2$-graded generator and vice
versa \cite[p.286]{L94}.

The flat $K$-group is related to other ordinary $K$-groups by the following exact
sequence \cite[\S 7.21]{Kar87}, \cite[(13)]{L94}
\begin{equation}\label{eq 2.2.1}
\begin{CD}
K^{-1}(X) @>r\circ\ch^{\odd}>> H^{\odd}(X; \R) @>\alpha>> K^{-1}_{\LL}(X) @>\beta>> K(X)
\end{CD}
\end{equation}
where $\ch^{\odd}$ is the odd Chern character, $r$ is induced by the inclusion of
coefficients $\Q\hookrightarrow\R$, and the maps $\alpha$ and $\beta$ are given by
\begin{displaymath}
\begin{split}
\alpha([\omega])&=(\C^n, h, \nabla^{\fl}, \omega)-(\C^n, h, \nabla^{\fl}, 0),\\
\beta(\E-\F)&=[E]-[F],
\end{split}
\end{displaymath}
where $\C^n\to X$ denotes the trivial complex vector bundle of rank $n$. As in the case
of ordinary $K$-theory, there exists a unique Chern character $\ch_{\R/\Q}:K^{-1}_{\LL}
(X)\to H^{\odd}(X; \R/\Q)$, called the flat Chern character \cite[Definition 9]{L94},
defined as follows. For a generator $\E=(E, h^E, \nabla^E, \omega)$ of $K^{-1}_{\LL}(X)$,
write $N=\rk(E)$. The condition $\ch(\nabla^E)-N=-d\omega$ implies $\ch(E-\C^N)=0\in
H^{\even}(X; \Q)$. Thus there exists $k\in\N$ such that $kE\cong k\C^N$. Let
$\nabla^{kE}_0$ be a unitary connection on $kE\to X$ with trivial holonomy. One can check
that the odd form \dis{\frac{1}{k}\CS(k\nabla, \nabla^{kE}_0)+\omega} is closed. The flat
Chern character $\ch_{\R/\Q}(\E)$ is defined to be
\begin{equation}\label{eq 2.2.2}
\ch_{\R/\Q}(\E)=\bigg[\frac{1}{k}\CS(k\nabla^E, \nabla^{kE}_0)+\omega\bigg]\mod\Q.
\end{equation}
Note that $\ch_{\R/\Q}(\E)$ is independent of the choices of $k$ and $\nabla^{kE}_0$
\cite[Lemma 1]{L94} and is a well defined group homomorphism \cite[Proposition 1]{L94}.

The flat Chern character of a $\Z_2$-graded generator $\E$ of the form (\ref{eq 2.2.0})
is defined as follows. The condition $\ch(\nabla^+)-\ch(\nabla^-)=-d\omega$ implies the
existence of $k\in\N$ such that $kE^+\cong kE^-$. Choose an isometric isomorphism
$j:kE^+\to kE^-$. Then $\ch_{\R/\Q}(\E)$ is defined to be
\begin{equation}\label{eq 2.2.3}
\ch_{\R/\Q}(\E)=\bigg[\frac{1}{k}\CS(k\nabla^+, j^*k\nabla^-)+\omega\bigg]\mod\Q,
\end{equation}
Note that $\ch_{\R/\Q}(\E)$ is independent of the choices of $k$ and $j$ \cite[p.289]{L94}.

\subsection{Local family index theorem}

In this subsection we recall the setup and the statement of the local FIT. We refer to
\cite{BGV} and the references therein for details.

Let $\pi:X\to B$ be a submersion with closed $\spin^c$ fibers of even relative dimension.
Denote by $T^VX\to X$ its vertical tangent bundle. Put a metric $g^{T^VX}$ on $T^VX\to X$.
Given a horizontal distribution $T^HX\to X$ and a Riemannian metric $g^{TB}$ on $TB\to B$,
we can define a metric on $TX\to X$ by $g^{TX}:=g^{T^VX}\oplus\pi^*g^{TB}$. If $\nabla^{TX}$
is the corresponding Levi-Civita connection, then $\nabla^{T^VX}:=P\circ\nabla^{TX}\circ P$
is a connection on $T^VX\to X$, where $P:TX\to T^VX$ is the orthogonal projection. Denote
by $S^c(T^VX)\to X$ the $\Z_2$-graded $\spin^c$ bundle and by $L^VX\to X$ the associated
characteristic Hermitian line bundle with a unitary connection $\nabla^{L^VX}$. Note that
the connection $\nabla^{T^VX}$ lifts uniquely to the local spinor bundle and preserves its
grading and the isomorphism $S^c(T^VX)\cong S(T^VX)\otimes L^VX$ exists globally
\cite[p.397]{LM89}. The connection $\nabla^{S^c(T^VX)}$ on $S^c(T^VX)\to X$, defined by
$$\nabla^{S^c(T^VX)}:=\nabla^{T^VX}\otimes\nabla^{L^VX},$$
preserves the grading of $S^c(T^VX)\to X$. The Todd form $\todd(\nabla^{S^c(T^VX)})$ of
$S^c(T^VX)\to X$ is defined to be
$$\todd(\nabla^{S^c(T^VX)})=\wh{A}(\nabla^{T^VX})\wedge e^{\frac{1}{2}c_1(\nabla^{L^V
X})}.$$

Define an infinite-rank bundle $\pi_*E\to B$ whose fiber over $z\in B$ is given by
$$(\pi_*E)_z:=\Gamma(X_z, (S^c(T^VX)\otimes E)_z).$$
Since $S^c(T^VX)\otimes E\to X$ is $\Z_2$-graded whose even and odd part are given by
\begin{equation}\label{eq 2.3.1}
(S^c(T^VX)\otimes E)^\pm=S^c(T^VX)^\pm\otimes E,
\end{equation}
it follows that the bundle $\pi_*E\to B$ is also $\Z_2$-graded whose even and odd part
are given by
\begin{equation}\label{eq 2.3.2}
(\pi_*E)^\pm_z=\Gamma(X_z, S^c(T^VX)^\pm\otimes E)_z)
\end{equation}
for each $z\in B$. The space of sections of $\pi_*E\to B$ is defined to be
\begin{equation}\label{eq 2.3.3}
\Gamma(B, \pi_*E):=\Gamma(X, S^c(T^VX)\otimes E).
\end{equation}
Note that $\pi_*E\to B$ admits an $L^2$-metric and a $\Z_2$-graded unitary connection
$\nabla^{\pi_*E}$ \cite[Proposition 9.13]{BGV}.

The $\spin^c$ Dirac operator $\DD^E:\Gamma(X, S^c(T^VX)\otimes E)\to\Gamma(X, S^c(T^VX)
\otimes E)$ is an odd operator given by
\begin{equation}\label{eq 2.3.4}
\DD^E=\sum_kc(e^k)\nabla_{e_k}^{S^c(T^VX)\otimes E},
\end{equation}
where $c$ is the Clifford multiplication, $\nabla^{S^c(T^VX)\otimes E}:=
\nabla^{S^c(T^VX)}\otimes\nabla^E$, $\set{e_k}$ is a local orthonormal frame for
$T^VX\to X$ and $\set{e^k}$ its dual frame for $(T^VX)^*\to X$. By (\ref{eq 2.3.3}),
$\DD^E$ can be regarded as an odd operator on $\pi_*E\to B$. Assume that  the family of
complex vector spaces $\ker(\DD^E_z)$ has locally constant dimension for $z\in B$. Then
$\ker(\DD^E_z)$ form a finite-rank $\Z_2$-graded complex vector bundle over $B$, denoted
by $\ker(\DD^E)\to B$ and is called the index bundle of $E\to X$. The analytic index
$\ind^{\an}(E)$ of $E\to X$ is defined by $\ind^{\an}(E)=[\ker(\DD^E)]\in K(B)$, and is
a ring homomorphism $\ind^{\an}:K(X)\to K(B)$.

Write $\EE=(E, h^E, \nabla^E)$. The Bismut superconnection $\bbb(\EE)$ on $\pi_*E\to B$
is defined to be
$$\bbb(\EE)=\DD^E+\nabla^{\pi_*E}-\frac{c(T)}{4},$$
where $T$ is the curvature 2-form of the fiber bundle $X\to B$. For each $z\in B$,
denote by $P_0^z:(\pi_*E)_z\to\ker(\DD^E)_z$ the orthogonal projection. Then $P_0$ is a
family of smoothing operators. Note that $\nabla^{\ker(\DD^E)}:=P_0\bbb(\EE)_{[1]}P_0$
is a $\Z_2$-graded unitary connection on $\ker(\DD^E)\to B$ \cite[Lemma 9.18]{BGV}. The
rescaled Bismut superconnection $\bbb(\EE)_t$ is defined to be
$$\bbb(\EE)_t=\sqrt{t}\DD^E+\nabla^{\pi_*E}-\frac{c(T)}{4\sqrt{t}}.$$
By \cite[Theorem 10.32]{BGV}, we have
\begin{eqnarray}
\lim_{t\to 0}\ch(\bbb(\EE)_t)&=&\int_{X/B}\todd(\nabla^{S^c(T^VX)})\wedge\ch
(\nabla^E), \label{eq 2.3.5}\\
\lim_{t\to\infty}\ch(\bbb(\EE)_t)&=&\ch(\nabla^{\ker(\DD^E)}).\label{eq 2.3.6}
\end{eqnarray}
Note that
\begin{equation}\label{eq 2.3.7}
\frac{d\ch(\bbb(\EE)_t)}{dt}=-d\str\bigg(\frac{d\bbb(\EE)_t}{dt}e^{-\frac{1}{2\pi i}
(\bbb(\EE)_t)^2}\bigg)
\end{equation}
and the integral \dis{\int^\infty_0\str\bigg(\frac{d\bbb(\EE)_t}{dt}e^{-\frac{1}
{2\pi i}(\bbb(\EE)_t)^2}\bigg)dt} converges \cite[Theorem 10.32]{BGV}. The eta form
\cite{BC89, D91} of $\EE$ is defined to be
\begin{equation}\label{eq 2.3.8}
\wt{\eta}(\EE):=\int^\infty_0\str\bigg(\frac{d\bbb(\EE)_t}{dt}e^{-\frac{1}{2\pi i}
(\bbb(\EE)_t)^2}\bigg)dt.
\end{equation}
The local FIT \cite[Theorem 10.32]{BGV} states that
\begin{equation}\label{eq 2.3.9}
d\wt{\eta}(\EE)=\int_{X/B}\todd(\nabla^{S^c(T^VX)})\wedge\ch(\nabla^E)-\ch
(\nabla^{\ker(\DD^E)}).
\end{equation}
which follows from (\ref{eq 2.3.5})--(\ref{eq 2.3.8}).
\begin{remark}
We use the slightly unconventional symbol $\bbb(\EE)$ and $\wt{\eta}(\EE)$ to emphasize
the dependence of the Bismut superconnection and the eta form on $\EE$. Of course they
also depend on other data: the metrics $g^{T^VX}$ and $g^{L^VX}$, the horizontal
distribution $T^HX$ and the unitary connection $\nabla^{L^VX}$. Henceforth we choose
and fix these data. Because of the definition of $K^{-1}_{\LL}$ we are only interested
in the deformation of the unitary connection on $E\to X$.
\end{remark}

\subsection{The flat analytic index}

In this subsection we recall the definition of the flat analytic index \cite[Definition
13]{L94}. Given a $\Z_2$-graded generator $\E$ of $K^{-1}_{\LL}(X)$, its flat analytic
index $\ind^{\an}_{\LL}(\E)\in K^{-1}_{\LL}(B)$ is, roughly speaking, given by the
analytic index of the $\Z_2$-graded data of $\E$ and a pushforward of the form $\omega$.
We refer to the construction of the analytic index in \S2.3, and indicate the changes
as follows.

Let $\pi:X\to B$ be a submersion with closed $\spin^c$ fibers of even relative
dimension, and $\E$ a $\Z_2$-graded generator of $K^{-1}_{\LL}(X)$ of the form
(\ref{eq 2.2.0}). As in \S2.3, the $\spin^c$ bundle $S^c(T^VX)\to X$ of $T^VX\to X$ is
$\Z_2$-graded. Since $E^+\oplus E^-\to X$ is also $\Z_2$-graded, the even and the odd
part of $S^c(T^VX)\wh{\otimes}E\to X$ become
$$(S^c(T^VX)\wh{\otimes}E)^\pm=S^c(T^VX)^+\otimes E^\pm\oplus S^c(T^VX)^-\otimes E^\mp.$$
It follows from (\ref{eq 2.3.2}) that the even and the odd part of $\pi_*E\to B$ has a
similar decomposition, so the same is true for $\ker(\DD^E)\to B$, that is,
$$\ker(\DD^E)^\pm=\ker(\DD^{E^+})^\pm\oplus\ker(\DD^{E^-})^\mp.$$
The $L^2$-metric and the Bismut superconnection $\bbb(\EE)$ on $\pi_*E\to B$ are
defined accordingly. The $\Z_2$-graded unitary connection on $\ker(\DD^E)\to B$ given
by
$$\nabla^{\ker(\DD^E)}=\nabla^{\ker(\DD^E)^+}\oplus\nabla^{\ker(\DD^E)^-},$$
where
$$\nabla^{\ker(\DD^E)^\pm}:=\nabla^{\ker(\DD^{E^+})^\pm}\oplus\nabla^{\ker
(\DD^{E^-})^\mp},$$
is a direct sum of connections. In this case the local FIT takes the form
\begin{equation}\label{eq 2.4.1}
\ch(\nabla^{\ker(\DD^E)})=\int_{X/B}\todd(\nabla^{S^c(T^VX)})\wedge\ch(\nabla^+\oplus
\nabla^-)-d\wt{\eta}(\EE^+\oplus\EE^-),
\end{equation}
where $\EE^+\oplus\EE^-=(E^+\oplus E^-, h^+\oplus h^-, \nabla^+\oplus\nabla^-)$.

The flat analytic index $\ind^{\an}_{\LL}:K^{-1}_{\LL}(X)\to K^{-1}_{\LL}(B)$ of a
$\Z_2$-graded generator $\E$ of the form (\ref{eq 2.2.0}) is defined by
\begin{equation}\label{eq 2.4.2}
\begin{split}
&~~~~\ind^{\an}_{\LL}(\E)\\
&:=\bigg(\ker(\DD^E), h^{\ker(\DD^E)}, \nabla^{\ker(\DD^E)}, \int_{X/B}\todd
(\nabla^{S^c(T^VX)})\wedge\omega+\wt{\eta}(\EE^+\oplus\EE^-)\bigg).
\end{split}
\end{equation}
It follows from (\ref{eq 2.3.9}) that $\ind^{\an}_{\LL}(\E)\in K^{-1}_{\LL}(B)$ is a
$\Z_2$-graded generator.

\section{Main results}

In this section we will prove the main results of this paper.

\subsection{Some properties of the eta form}

In this subsection we provide proofs of the additivity and the variational formula of
the eta forms.

Although $\ind^{\an}_{\LL}$ is defined using the spin$^c$ Dirac operator, we prove the
additivity of the eta forms in a slightly more general setting. Instead of working on
the twisted $\spin^c$ bundle $S^c(T^VX)\otimes E\to X$ we work on Clifford modules.
Before we state and prove the result we briefly recall the Clifford modules in our
setup. We refer to \cite[\S10.2, \S10.3]{BGV} for the details.

Let $\pi:X\to B$ be a submersion with closed fibers of even relative dimension. Put a
Riemannian metric $g^{TB}$ on $TB\to B$ and a metric $g^{T^VX}$ on the vertical bundle
$T^VX\to X$. Recall from \cite[p.322]{BGV} that a Clifford module along the fibers of
$\pi:X\to B$ is given by $\e=(E, h^E, \nabla^E)$, where $E\to X$ is a $\Z_2$-graded
complex vector bundle, $h^E$ a $\Z_2$-graded Hermitian metric and $\nabla^E$ a
$\Z_2$-graded unitary connection, with a skew-adjoint action
$$c:\Cl((T^VX)^*)\to\End(E),$$
where $(T^VX)^*\to X$ denotes the dual bundle of $T^VX\to X$ and $\Cl((T^VX)^*)\to X$
is the Clifford bundle of $(T^VX)^*\to X$, such that
$$[\nabla^E_V, c(\alpha)]=c(\nabla^{(T^VX)^*}_V\alpha)$$
for $V\in\Gamma(X, TX)$ and $\alpha\in\Gamma(X, (T^VX)^*)$. One can define a Dirac
operator $\DD^\e:\Gamma(X, \e)\to\Gamma(X, \e)$ for a Clifford module $\e$ along the
fibers in a way similar to (\ref{eq 2.3.4}).

Here we recall the definition of the Bismut superconnection $\bbb(\e)$ associated to
$\e$ \cite[Proposition 10.15]{BGV}. Given a Clifford module $\e=(E, h^E, \nabla^E)$
along the fibers of $X\to B$, define a complex vector bundle $\wt{E}:=\pi^*\Lambda
(T^*B)\otimes E$ over $X$ and equip it with the Hermitian metric $\pi^*g^{T^*B}\otimes
h^E$. Consider the Clifford algebra bundle $\cl_0(T^*X)\to X$, where $\cl_0(T^*X)$
denotes $\cl(T^*X)$ equipped with the degenerate metric $g^{T^*X}_0$, with
$g^{T^*X}_\varepsilon:=g^{(T^VX)^*}\oplus\varepsilon\pi^*g^{T^*B}$. Since $T^*X\cong
(T^VX)^*\oplus\pi^*T^*B$, it follows that $\cl_0(T^*X)\cong\pi^*\Lambda(T^*B)\otimes
\cl((T^VX)^*)$. Define a Clifford multiplication $m_0:\cl_0(T^*X)\to\End(\wt{E})$ by
$m_0(\alpha)=\alpha\wedge\cdot$ if $\alpha\in\Gamma(X, \pi^*T^*B)$ and $m_0(\alpha)=
c(\alpha)$ if $\alpha\in\Gamma(X, (T^VX)^*)$. Define a connection $\nabla^{\wt{E}}$
on $\wt{E}\to X$ by
\begin{equation}\label{eq 3.1.1}
\nabla^{\wt{E}}=\pi^*\nabla^{T^*B}\otimes\nabla^E+\frac{1}{2}m_0(\omega),
\end{equation}
where $\omega\in\Omega^1(X, \Lambda^2(T^*X))$ is characterized by
\cite[Proposition 10.6]{BGV}. By \cite[Proposition 10.10]{BGV} $\nabla^{\wt{E}}$ is a
Clifford connection. Thus $\wt{\e}:=(\wt{E}, \pi^*g^{T^*B}\otimes h^E, \nabla^{\wt{E}})$
is a Clifford module over the Clifford algebra bundle $\cl_0(T^*X)\to X$. Note that
$\Omega(B, \pi_*E)$ is defined to be $\Gamma(X, \wt{E})$. The Bismut superconnection
$\bbb(\e):\Omega(B, \pi_*E)\to\Omega(B, \pi_*E)$ is defined as a Dirac operator $\bbb
(\e):\Gamma(X, \wt{E})\to\Gamma(X, \wt{E})$ by the formula
\begin{equation}\label{eq 3.1.2}
\bbb(\e)=\sum_km_0(e^k)\nabla^{\wt{E}}_{e_k},
\end{equation}
where $\set{e_k}$ is a local orthonormal frame for $TX\to X$ and $\set{e^k}$ its dual
frame for $T^*X\to X$.
\begin{prop}\label{prop 1}
Let $\pi:X\to B$ be a submersion with closed fibers of even relative dimension,
$\e=(E, h^E, \nabla^E)$ a Clifford module over $X$ along the fibers of $X\to B$ and
$\DD^\e$ the Dirac operator associated to $\e$. Let $\wt{\eta}(\e)$ be the eta form
of the Bismut superconnection $\bbb(\e)$. If $\f=(F, h^F, \nabla^F)$ is another
Clifford module over $X$ along the fibers of $X\to B$, then
$$\wt{\eta}(\e\oplus\f)=\wt{\eta}(\e)+\wt{\eta}(\f)$$
up to exact forms.
\end{prop}
\begin{proof}
First of all we claim that $\bbb(\e\oplus\f)=\bbb(\e)\oplus\bbb(\f)$. By (\ref{eq 3.1.2})
it suffices to prove that
$$\nabla^{\wt{E\oplus F}}=\nabla^{\wt{E}}\oplus\nabla^{\wt{F}}.$$
To see this, let $\beta\otimes(\alpha_1\oplus\alpha_2)\in\Gamma(X, \pi^*\Lambda(T^*B)
\otimes(E\oplus F))$. By (\ref{eq 3.1.1}) we have
\begin{displaymath}
\begin{split}
&\qquad\nabla^{\wt{E\oplus F}}(\beta\otimes(\alpha_1\oplus \alpha_2))\\
&=(\pi^*\nabla^{TB}\otimes\nabla^{E\oplus F})(\beta\otimes(\alpha_1\oplus\alpha_2))+
\frac{1}{2}m_0(\omega)(\beta\otimes(\alpha_1\oplus \alpha_2))\\
&=\pi^*\nabla^{TB}\beta\otimes(\alpha_1\oplus \alpha_2)+\beta\otimes(\nabla^E\alpha_1
\oplus\nabla^F\alpha_2)+\frac{1}{2}m_0(\omega)\beta\otimes(\alpha_1\oplus\alpha_2)\\
&=\bigg(\pi^*\nabla^{TB}\beta\otimes\alpha_1+\beta\otimes\nabla^E\alpha_1+\frac{1}{2}m_0
(\omega)\beta\otimes\alpha_1\bigg)\\
&\qquad\oplus\bigg(\pi^*\nabla^{TB}\beta\otimes\alpha_2+\beta\otimes\nabla^F\alpha_2+
\frac{1}{2}m_0(\omega)\beta\otimes\alpha_2\bigg)\\
&=\nabla^{\wt{E}}(\beta\otimes\alpha_1)\oplus\nabla^{\wt{F}}(\beta\otimes\alpha_2)=
(\nabla^{\wt{E}}\oplus\nabla^{\wt{F}})(\beta\otimes(\alpha_1\oplus\alpha_2)).
\end{split}
\end{displaymath}
The additivity of the Bismut superconnections holds for the rescaled Bismut
superconnection; i.e.,
$$\bbb(\e\oplus\f)_t=\bbb(\e)_t\oplus\bbb(\f)_t.$$
Consider the Chern--Simons form
$$\CS(\bbb(\e)_T, \bbb(\e)_t)=\int^T_t\str\bigg(\frac{d\bbb(\e)_s}{ds}e^{-\frac{1}
{2\pi i}(\bbb(\e)_s)^2}\bigg)ds,$$
where $0<t<T$ are fixed. Properties (\ref{eq 2.1.6})-(\ref{eq 2.1.8}) extend to this
case. Therefore
\begin{displaymath}
\begin{split}
\CS(\bbb(\e\oplus\f)_T, \bbb(\e\oplus\f)_t)&=\CS(\bbb(\e)_T\oplus\bbb(\f)_T,
\bbb(\e)_t\oplus\bbb(\f)_t)\\
&=\CS(\bbb(\e)_T, \bbb(\e)_t)+\CS(\bbb(\f)_T, \bbb(\f)_t).
\end{split}
\end{displaymath}
By letting $T\to\infty$ and $t\to 0$ in above, the convergence of all the integrals
involved \cite[Theorem 10.32]{BGV} shows that $\wt{\eta}(\e\oplus\f)=\wt{\eta}(\e)+
\wt{\eta}(\f)$ up to exact forms.
\end{proof}
Let $\EE=(E, h^E, \nabla^E)$ and $\FF=(F, h^F, \nabla^F)$ be Hermitian bundles with
unitary connections. By applying Proposition \ref{prop 1} to the twisted $\spin^c$
bundle $S^c(T^VX)\otimes(E\oplus F)\to X$ where the fibers of $\pi:X\to B$ are assumed
to be spin$^c$, we have
\begin{equation}\label{eq 3.1.3}
\wt{\eta}(\EE\oplus\FF)=\wt{\eta}(\EE)+\wt{\eta}(\FF)
\end{equation}
up to exact forms.

Some remarks for Proposition \ref{prop 2}.
\begin{remark}\label{remark 2}
Let $p:M\to X$ be a smooth fiber bundle with compact fibers. By \cite[Chapter 1]{BT82},
we have
\begin{equation}\label{eq 3.1.4}
\int_{M/X}p^*\alpha\wedge\beta=\alpha\wedge\bigg(\int_{M/X}\beta\bigg),
\end{equation}
for all $\alpha\in\Omega(X)$ and $\beta\in\Omega(M)$.

If $q:X\to B$ is another smooth fiber bundle with compact fibers, then $q\circ p:M\to B$
is a smooth fiber bundle with compact fibers, then it is straightforward to check (or
see \cite[Problem 3 (p.311)]{GHV}) that
\begin{equation}\label{eq 3.1.5}
\int_{M/B}=\int_{X/B}\circ\int_{M/X}.
\end{equation}
\end{remark}
\begin{prop}\label{prop 2}
Let $\pi:X\to B$ be a submersion with closed $\spin^c$ fibers of even relative dimension.
Write $\EE_k=(E, h^E, \nabla_k^E)$, where $k\in\set{0, 1}$, as in \S2.1. Then
\begin{equation}\label{eq 3.1.6}
\wt{\eta}(\EE_1)-\wt{\eta}(\EE_0)=\int_{X/B}\todd(\nabla^{S^c(T^VX)})\wedge\CS(\nabla^E_1,
\nabla^E_0)-\CS(\nabla^{\ker(\DD^E)}_1, \nabla^{\ker(\DD^E)}_0)
\end{equation}
up to exact forms.
\end{prop}
Proposition \ref{prop 2} is a special case of the variational formula of the equivariant
eta forms \cite[Theorem 1.7]{L15}, where the geometric data on $\pi:X\to B$ and the
connection on $E\to X$ are deformed.
\begin{proof}
Consider the following commutative diagram
\cdd{\wt{X} @>p_X>> X \\ @V\wt{\pi}VV @VV\pi V \\ \wt{B} @>>p_B> B}
The geometric data on $\wt{\pi}:\wt{X}\to\wt{B}$ is obtained by pulling back the
geometric data on $\pi:X\to B$. Then the local FIT (\ref{eq 2.3.9})
for $\wt{\EE}=(\wt{E}, h^{\wt{E}}, \nabla^{\wt{E}})$ gives
\begin{equation}\label{eq 3.1.7}
d\wt{\eta}(\wt{\EE})=\int_{\wt{X}/\wt{B}}\todd(\nabla^{S^c(T^V\wt{X})})\wedge\ch
(\nabla^{\wt{E}})-\ch(\nabla^{\ker(\DD^{\wt{E}})}).
\end{equation}
Consider $\ker(\DD^{\wt{E}})\to\wt{B}$. By the same reason as in \S2.1 we have
$i_0^*\ker(\DD^{\wt{E}})\cong i_1^*\ker(\DD^{\wt{E}})$. Moreover, for $k\in\set{0, 1}$
the connection defining the $\spin^c$ Dirac operator $\DD^{\wt{E}|_{X\times\set{k}}}$
on $(S^c(T^V\wt{X})\otimes\wt{E})|_{X\times\set{k}}\to X\times\set{k}$ is $\nabla^{S^c
(T^V\wt{X})}\otimes\nabla_k$. Thus $\ker(\DD^{E_k})\cong i_k^*\ker(\DD^{\wt{E}})$, so
$\ker(\DD^{E_0})\cong\ker(\DD^{E_1})$, and is therefore denoted by $\ker(\DD^E)$.
Write $\nabla^{\ker(\DD^E)}_k$ for the unitary connection on $\ker(\DD^E)\to B$ induced
by $\nabla_k$. Denote by $i:\partial\wt{B}\to\wt{B}$ the inclusion map. By
(\ref{eq 2.1.3}), we have
$$\wt{\eta}(\EE_1)-\wt{\eta}(\EE_0)=\int_{\partial\wt{B}/B}i^*\wt{\eta}(\wt{\EE})
=\int_{\wt{B}/B}d_{\wt{B}}\wt{\eta}(\wt{\EE})-d_B\int_{\wt{B}/B}\wt{\eta}(\wt{\EE}).$$
By modding out exact forms, it follows from (\ref{eq 3.1.7}) that
\begin{displaymath}
\begin{split}
\wt{\eta}(\EE_1)-\wt{\eta}(\EE_0)&=\int_{\wt{B}/B}d_{\wt{B}}\wt{\eta}(\wt{\EE})\\
&=\int_{\wt{B}/B}\bigg(\int_{\wt{X}/\wt{B}}\todd(\nabla^{S^c(T^V\wt{X})})\wedge\ch
(\nabla^{\wt{E}})-\ch(\nabla^{\ker(\DD^{\wt{E}})})\bigg)\\
&=\int_{\wt{B}/B}\int_{\wt{X}/\wt{B}}p_X^*\todd(\nabla^{S^c(T^VX)})\wedge\ch
(\nabla^{\wt{E}})-\int_{\wt{B}/B}\ch(\nabla^{\ker(\DD^{\wt{E}})}).
\end{split}
\end{displaymath}
By (\ref{eq 2.1.2}), the last term of the right-hand side is equal to
$\CS(\nabla_1^{\ker(\DD^E)}, \nabla_0^{\ker(\DD^E)})$. Then
\begin{displaymath}
\begin{split}
&\qquad\wt{\eta}(\EE_1)-\wt{\eta}(\EE_0)\\
&=\int_{\wt{B}/B}\int_{\wt{X}/\wt{B}}p_X^*\todd(\nabla^{S^c(T^VX)})\wedge\ch
(\nabla^{\wt{E}})-\CS(\nabla_1^{\ker(\DD^E)}, \nabla_0^{\ker(\DD^E)})\\
&=\int_{\wt{X}/B}p_X^*\todd(\nabla^{S^c(T^VX)})\wedge\ch(\nabla^{\wt{E}})-
\CS(\nabla_1^{\ker(\DD^E)}, \nabla_0^{\ker(\DD^E)})\\
&=\int_{X/B}\int_{\wt{X}/X}p_X^*\todd(\nabla^{S^c(T^VX)})\wedge\ch(\nabla^{\wt{E}})
-\CS(\nabla_1^{\ker(\DD^E)}, \nabla_0^{\ker(\DD^E)})\\
&=\int_{X/B}\todd(\nabla^{S^c(T^VX)})\wedge\bigg(\int_{\wt{X}/X}\ch(\nabla^{\wt{E}})
\bigg)-\CS(\nabla_1^{\ker(\DD^E)}, \nabla_0^{\ker(\DD^E)})\\
&=\int_{X/B}\todd(\nabla^{S^c(T^VX)})\wedge\CS(\nabla^E_1, \nabla^E_0)-\CS(\nabla_1^{\ker
(\DD^E)}, \nabla_0^{\ker(\DD^E)}),
\end{split}
\end{displaymath}
up to exact forms, where the second and the third equalities follow from (\ref{eq 3.1.5}),
the forth equality follows from (\ref{eq 3.1.4}) and the last equality follows from
(\ref{eq 2.1.2}).
\end{proof}
We call (\ref{eq 3.1.6}) the variational formula of the eta forms of the pair $(\EE_1,
\EE_0)$.

\subsection{The flat GRR}

In this subsection we prove that the flat analytic index $\ind^{\an}_{\LL}:K^{-1}_{\LL}(X)
\to K^{-1}_{\LL}(B)$ is well defined and the flat GRR. The proof of Proposition
\ref{prop 3} is the essentially the same as \cite[Proposition 3]{H14}.
\begin{prop}\label{prop 3}
Let $\pi:X\to B$ be a submersion with closed $\spin^c$ fibers of even relative dimension.
The flat analytic index
$$\ind^{\an}_{\LL}:K^{-1}_{\LL}(X)\to K^{-1}_{\LL}(B)$$
is well defined.
\end{prop}
\begin{proof}
For $k=0, 1$, let $\E_k=(E_k^+\oplus E_k^-, h^+_k\oplus h^-_k, \nabla^+_k\oplus\nabla^-_k,
\omega_k)$ be $\Z_2$-graded generators of $K^{-1}_{\LL}(X)$ such that the classes of $\E_0$
and $\E_1$ are equal in $K^{-1}_{\LL}(X)$. For notational clarity we write $E_k$ for
$E_k^+\oplus E_k^-$ and similarly for other $\Z_2$-graded objects. By the definition of
$K^{-1}_{\LL}$, there exists a $\Z_2$-graded generator $\G=(G, h^G, \nabla^G, \omega^G)$
of $K^{-1}_{\LL}(X)$ such that
\begin{equation}\label{eq 3.2.1}
E_1\oplus G\cong E_0\oplus G,
\end{equation}
and
\begin{equation}\label{eq 3.2.2}
\omega_1-\omega_0=\CS(\nabla_0\oplus\nabla^G, \nabla_1\oplus\nabla^G).
\end{equation}
Since the analytic index is additive, (\ref{eq 3.2.1}) implies
\begin{equation}\label{eq 3.2.3}
\begin{split}
\ind^{\an}(E_1)\oplus\ind^{\an}(G)&\cong\ind^{\an}(E_1\oplus G)\\
&\cong\ind^{\an}(E_0\oplus G)\cong\ind^{\an}(E_0)\oplus\ind^{\an}(G).
\end{split}
\end{equation}
Since the diagram
\cdd{K^{-1}(X) @>\beta>> K(X) \\ @V\ind^{\an}_{\LL}VV @VV\ind^{\an}V \\ K^{-1}(B)
@>>\beta> K(B)}
commutes, it follows that
$$\beta(\ind^{\an}_{\LL}(\E_1)-\ind^{\an}_{\LL}(\E_0))=\ind^{\an}(E_1)-\ind^{\an}
(E_0)=0.$$
It follows from the exact sequence (\ref{eq 2.2.1}) that there exists $[\omega]\in
H^{\odd}(B; \R)$ such that
\begin{equation}\label{eq 3.2.4}
\alpha([\omega])=\ind^{\an}_{\LL}(\E_1)-\ind^{\an}_{\LL}(\E_0).
\end{equation}
It suffices to prove that $\omega$ is an exact form. By the definition of $\alpha$ and
(\ref{eq 2.4.2}), (\ref{eq 3.2.4}) implies
\begin{displaymath}
\begin{split}
&\bigg(\ker(\DD^{E_1}), h^{\ker(\DD^{E_1})}, \nabla^{\ker(\DD^{E_1})}, \int_{X/B}
\todd(\nabla^{S^c(T^VX)})\wedge\omega_1+\wt{\eta}(\EE_1)\bigg)\\
&\qquad-\bigg(\ker(\DD^{E_0}), h^{\ker(\DD^{E_0})}, \nabla^{\ker(\DD^{E_0})}, \int_{X/B}
\todd(\nabla^{S^c(T^VX)})\wedge\omega_0+\wt{\eta}(\EE_0)\bigg)\\
&=(\C^n, h, \nabla^{\fl}, \omega)-(\C^n, h, \nabla^{\fl}, 0).
\end{split}
\end{displaymath}
It follows from the definition of $K^{-1}_{\LL}$ that
\begin{equation}\label{eq 3.2.5}
\begin{split}
\omega&=\CS(\nabla^{\ker(\DD^{E_1})}\oplus\nabla^{\fl}, \nabla^{\ker(\DD^{E_0})}\oplus
\nabla^{\fl})+\wt{\eta}(\EE_1)-\wt{\eta}(\EE_0)\\
&~~~~~~~~+\int_{X/B}\todd(\nabla^{S^c(T^VX)})\wedge(\omega_1-\omega_0)
\end{split}
\end{equation}
up to exact forms. By the additivity of the eta forms (\ref{eq 3.1.3}), we have
$$\wt{\eta}(\EE_1)-\wt{\eta}(\EE_0)=\wt{\eta}(\EE_1\oplus\GG)-\wt{\eta}(\EE_0\oplus\GG).$$
Together with (\ref{eq 2.1.8}) the right-hand side of (\ref{eq 3.2.5}) becomes
\begin{displaymath}
\begin{split}
&\CS(\nabla^{\ker(\DD^{E_1})}\oplus\nabla^{\ker(\DD^G)}\oplus\nabla^{\fl}, \nabla^{\ker
(\DD^{E_0})}\oplus\nabla^{\ker(\DD^G)}\oplus\nabla^{\fl})\\
&~~~~+\int_{X/B}\todd(\nabla^{S^c(T^VX)})\wedge\CS(\nabla_0\oplus\nabla^G, \nabla_1
\oplus\nabla^G)+\wt{\eta}(\EE_1\oplus\GG)-\wt{\eta}(\EE_0\oplus\GG).
\end{split}
\end{displaymath}
By (\ref{eq 3.2.1}) and (\ref{eq 2.1.8}), it becomes
\begin{displaymath}
\begin{split}
&\CS(\nabla^{\ker(\DD^{E_1})}\oplus\nabla^{\ker(\DD^G)}, \nabla^{\ker(\DD^{E_0})}\oplus
\nabla^{\ker(\DD^G)})+\wt{\eta}(\EE_1\oplus\GG)-\wt{\eta}(\EE_0\oplus\GG)\\
&~~~~+\int_{X/B}\todd(\nabla^{S^c(T^VX)})\wedge\CS(\nabla_0\oplus\nabla^G, \nabla_1\oplus
\nabla^G).
\end{split}
\end{displaymath}
Because of (\ref{eq 3.2.3}), we can apply the variational formula for the eta forms
(Proposition \ref{prop 2}) to the pair $(\EE_0\oplus\GG, \EE_1\oplus\GG)$, which shows
that the above form is exact. From (\ref{eq 3.2.5}) we see that $\omega$ is exact, so
$\ind^{\an}_{\LL}(\E_1)=\ind^{\an}_{\LL}(\E_0)$. Therefore the flat analytic index
$\ind^{\an}_{\LL}$ is well defined.
\end{proof}

We are now ready to prove the flat GRR.
\begin{thm}\label{thm 1}
Let $\pi:X\to B$ be a submersion with closed $\spin^c$ fibers of even relative dimension.
The following diagram commutes.
\cdd{K^{-1}_{\LL}(X) @>\ch_{\R/\Q}>> H^{\odd}(X; \R/\Q) \\ @V\ind^{\an}_{\LL}VV
@VV\int_{X/B}\todd(X/B)\cup(\cdot) V \\ K^{-1}_{\LL}(B) @>>\ch_{\R/\Q}>
H^{\odd}(B; \R/\Q)}
i.e., for a $\Z_2$-graded generator $\E$ of $K^{-1}_{\LL}(X)$ of the form
(\ref{eq 2.2.0}), we have
\begin{equation}\label{eq 3.2.6}
\ch_{\R/\Q}(\ind^{\an}_{\LL}(\E))=\int_{X/B}\todd(X/B)\cup\ch_{\R/\Q}(\E).
\end{equation}
\end{thm}
\begin{proof}
By (\ref{eq 2.2.3}) and (\ref{eq 2.4.2}), $\ch_{\R/\Q}(\ind^{\an}_{\LL}(\E))$ is given by
the mod $\Q$ reduction of the de Rham class of
\begin{equation}\label{eq 3.2.7}
\frac{1}{\ell}\CS(\ell\nabla^{\ker(\DD^E)^+}, j_1^*\ell\nabla^{\ker(\DD^E)^-})+
\int_{X/B}\todd(\nabla^{S^c(T^VX)})\wedge\omega+\wt{\eta}(\EE),
\end{equation}
where $\ell\in\N$ and $j_1:\ell\ker(\DD^E)^+\to\ell\ker(\DD^E)^-$ is an isometric
isomorphism. Similarly \dis{\int_{X/B}\todd(X/B)\cup\ch_{\R/\Q}(\E)} is given by the mod
$\Q$ reduction of the de Rham class of
\begin{equation}\label{eq 3.2.8}
\frac{1}{k}\int_{X/B}\todd(\nabla^{S^c(T^VX)})\wedge\CS(k\nabla^+, j_2^*k\nabla^-)+
\int_{X/B}\todd(\nabla^{S^c(T^VX)})\wedge\omega,
\end{equation}
where $k\in\N$ and $j_2:kE^+\to kE^-$ is an isometric isomorphism. Consider the
difference between (\ref{eq 3.2.7}) and (\ref{eq 3.2.8}), which is given by
\begin{equation}\label{eq 3.2.9}
\begin{split}
h:&=\frac{1}{\ell}\CS(\ell\nabla^{\ker(\DD^E)^+}, j_1^*\ell\nabla^{\ker(\DD^E)^-})+
\wt{\eta}(\EE)\\
&\qquad-\frac{1}{k}\int_{X/B}\todd(\nabla^{S^c(T^VX)})\wedge\CS(k\nabla^+, j_2^*k
\nabla^-).
\end{split}
\end{equation}
Thus to prove (\ref{eq 3.2.6}) it suffices to prove $h=0$ up to forms with periods in
$\Q$. Let $m$ be the least common multiple of $k$ and $\ell$. Then there exist unique
$d_1, d_2\in\N$ such that $m=\ell d_1$ and $m=kd_2$. Since $j_2:kE^+\to kE^-$ is an
isometric isomorphism, the same is true for
$$d_2j_2:\overbrace{kE^+\oplus\cdots\oplus kE^+}^{d_2}\to\overbrace{kE^-\oplus\cdots
\oplus kE^-}^{d_2},$$
where $d_2j_2:=\overbrace{j_2\oplus\cdots\oplus j_2}^{d_2}$. Moreover, $d_2(k\nabla^+)
=m\nabla^+$ and $(d_2j_2)^*(d_2k\nabla^+)=(d_2j_2)^*m\nabla^-$ are unitary connections
on $mE^+\to X$. Note that
\begin{equation}\label{eq 3.2.10}
\begin{split}
\frac{1}{k}\CS(k\nabla^+, j_2^*k\nabla^-)&=\frac{1}{m}\CS(m\nabla^+, (d_2j_2)^*m
\nabla^-)+\frac{1}{k}\CS(k\nabla^+, j_2^*k\nabla^-)\\
&\qquad-\frac{1}{m}\CS(m\nabla^+, (d_2j_2)^*m\nabla^-).\\
\end{split}
\end{equation}
By (\ref{eq 2.1.6}) and (\ref{eq 2.1.7}), the last two terms of the right-hand side
of (\ref{eq 3.2.10}) equal
\begin{displaymath}
\begin{split}
&\qquad\frac{1}{k}\CS(k\nabla^+, j_2^*k\nabla^-)-\frac{1}{m}\CS(m\nabla^+, (d_2j_2)^*m
\nabla^-)\\
&=\frac{1}{d_2k}\CS(d_2k\nabla^+, (d_2j_2)^*(d_2k\nabla^-))-\frac{1}{m}\CS(m\nabla^+,
(d_2j_2)^*m\nabla^-)\\
&=\frac{1}{m}\bigg(\CS(m\nabla^+, (d_2j_2)^*(m\nabla^-))-\CS(m\nabla^+, (d_2j_2)^*m
\nabla^-)\bigg)\\
&=\frac{1}{m}\bigg(\CS(m\nabla^+, (d_2j_2)^*(m\nabla^-))+\CS((d_2j_2)^*m\nabla^-,
m\nabla^+)\bigg)\\
&=\frac{1}{m}\CS(m\nabla^+, m\nabla^+)=0
\end{split}
\end{displaymath}
up to exact forms.\footnote{The argument is similar to the proof of \cite[Lemma 1]{L94},
which says that the definition of $\ch_{\R/\Q}$ is, in particular, independent of the
choice of $k\in\N$ (see (\ref{eq 2.2.2})).} Thus
\begin{displaymath}
\begin{split}
&\frac{1}{k}\int_{X/B}\todd(\nabla^{S^c(T^VX)})\wedge\CS(k\nabla^+, j_2^*k\nabla^-)\\
&\qquad=\frac{1}{m}\int_{X/B}\todd(\nabla^{S^c(T^VX)})\wedge\CS(m\nabla^+, (d_2j_2)^*
m\nabla^-)
\end{split}
\end{displaymath}
up to exact forms. By the same argument, the first term of the right-hand side of
(\ref{eq 3.2.9}) becomes
$$\frac{1}{\ell}\CS(\ell\nabla^{\ker(\DD^E)^+}, j_1^*\ell\nabla^{\ker(\DD^E)^-})=
\frac{1}{m}\CS(m\nabla^{\ker(\DD^E)^+}, (d_1j_1)^*m\nabla^{\ker(\DD^E)^-})$$
up to exact forms. Since exact forms have zero periods, it follows that (\ref{eq 3.2.9})
becomes
\begin{displaymath}
\begin{split}
h&=\frac{1}{m}\CS(m\nabla^{\ker(\DD^E)^+}, (d_1j_1)^*m\nabla^{\ker(\DD^E)^-})+\wt{\eta}
(\EE)\\
&\qquad-\frac{1}{m}\int_{X/B}\todd(\nabla^{S^c(T^VX)})\wedge\CS(m\nabla^+, (d_2j_2)^*
m\nabla^-).
\end{split}
\end{displaymath}
Since $\wt{\eta}(m\EE)=m\wt{\eta}(\EE)$ by (\ref{eq 3.1.3}), it follows that
\begin{displaymath}
\begin{split}
h&=\frac{1}{m}\bigg(\CS(m\nabla^{\ker(\DD^E)^+}, (d_1j_1)^*m\nabla^{\ker(\DD^E)^-})+
\wt{\eta}(m\EE)\\
&\qquad-\int_{X/B}\todd(\nabla^{S^c(T^VX)})\wedge\CS(m\nabla^+, (d_2j_2)^*m\nabla^-)
\bigg).
\end{split}
\end{displaymath}
\emph{A priori} the variational formula of the eta forms (Proposition \ref{prop 2})
cannot be applied to the pair $(m\EE^+, (d_2j_2)^*m\EE^-)$ since the isometric
isomorphisms $d_1j_1$ and $d_2j_2$ are not related in general. However, as remarked in
\cite[p.289]{L94}, the flat Chern character $\ch_{\R/\Q}$ of $\Z_2$-graded generator is
independent of the choice of the isometric isomorphism involved. Thus, without loss of
generality, we can assume $d_1j_1$ is induced by $d_2j_2$, so Proposition \ref{prop 2}
can be applied. Thus $h=0$ up to exact forms, and therefore (\ref{eq 3.2.6}) holds.

We have proved (\ref{eq 3.2.6}) under the assumption that the kernel bundle $\ker(\DD^E)
\to B$ exists. Without this assumption one can proceed as in \cite[\S5]{L94} and
\cite[\S7.12]{FL10} to prove (\ref{eq 3.2.6}).
\end{proof}
\bibliographystyle{amsplain}
\bibliography{MBib}

\providecommand{\bysame}{\leavevmode\hbox to3em{\hrulefill}\thinspace}
\providecommand{\MR}{\relax\ifhmode\unskip\space\fi MR }
\providecommand{\MRhref}[2]{%
  \href{http://www.ams.org/mathscinet-getitem?mr=#1}{#2}
}
\providecommand{\href}[2]{#2}
\begin{thebibliography}{10}

\bibitem{APS75b}
M.~F. Atiyah, V.~K. Patodi, and I.~M. Singer, \emph{Spectral asymmetry and
  {R}iemannian geometry \textrm{II}}, Math. Proc. Cambridge Philos. Soc.
  \textbf{78} (1975), 405--432.

\bibitem{APS76}
\bysame, \emph{Spectral asymmetry and {R}iemannian geometry. {III}}, Math.
  Proc. Cambridge Philos. Soc. \textbf{79} (1976), no.~1, 71--99.

\bibitem{AS71}
M.~F. Atiyah and I.~M. Singer, \emph{The index of elliptic operators. {IV}},
  Ann. of Math. (2) \textbf{93} (1971), 119--138.

\bibitem{AS84}
\bysame, \emph{Dirac operators coupled to vector potentials}, Proc. Nat. Acad.
  Sci. U.S.A. \textbf{81} (1984), no.~8, Phys. Sci., 2597--2600.

\bibitem{BD82}
Paul Baum and Ronald~G. Douglas, \emph{{$K$} homology and index theory},
  Operator algebras and applications, {P}art {I} ({K}ingston, {O}nt., 1980),
  Proc. Sympos. Pure Math., vol.~38, Amer. Math. Soc., Providence, R.I., 1982,
  pp.~117--173.

\bibitem{BGV}
Nicole Berline, Ezra Getzler, and Mich{\`e}le Vergne, \emph{Heat kernels and
  {D}irac operators}, Grundlehren Text Editions, Springer-Verlag, Berlin, 2004,
  Corrected reprint of the 1992 original.

\bibitem{B05}
Jean-Michel Bismut, \emph{Eta invariants, differential characters and flat
  vector bundles}, Chinese Ann. Math. Ser. B \textbf{26} (2005), 15--44.

\bibitem{BC89}
Jean-Michel Bismut and Jeff Cheeger, \emph{$\eta$-invariants and their
  adiabatic limits}, J. Amer. Math. Soc. \textbf{2} (1989), 33--70.

\bibitem{BF86a}
Jean-Michel Bismut and Daniel~S. Freed, \emph{The analysis of elliptic
  families. {I}. {M}etrics and connections on determinant bundles}, Comm. Math.
  Phys. \textbf{106} (1986), no.~1, 159--176.

\bibitem{BF86b}
\bysame, \emph{The analysis of elliptic families. {II}. {D}irac operators, eta
  invariants, and the holonomy theorem}, Comm. Math. Phys. \textbf{107} (1986),
  no.~1, 103--163.

\bibitem{BT82}
Raoul Bott and Loring~W. Tu, \emph{Differential forms in algebraic topology},
  Graduate Texts in Mathematics, vol.~82, Springer-Verlag, New York-Berlin,
  1982.

\bibitem{BS09}
Ulrich Bunke and Thomas Schick, \emph{Smooth {$K$}-theory}, Ast\'erisque
  (2009), no.~328, 45--135 (2010).

\bibitem{BS10}
\bysame, \emph{Uniqueness of smooth extensions of generalized cohomology
  theories}, J. Topol. \textbf{3} (2010), 110--156.

\bibitem{BS12}
\bysame, \emph{Differential {K}-theory: a survey}, Global differential
  geometry, Springer Proc. Math., vol.~17, Springer, Heidelberg, 2012,
  pp.~303--357.

\bibitem{C87}
Jeff Cheeger, \emph{{$\eta$}-invariants, the adiabatic approximation and
  conical singularities. {I}. {T}he adiabatic approximation}, J. Differential
  Geom. \textbf{26} (1987), no.~1, 175--221.

\bibitem{CS85}
Jeff Cheeger and James Simons, \emph{Differential characters and geometric
  invariants}, in Geometry and Topology (College Park, Md., 1983/84), Lecture
  Notes in Math. \textbf{1167} (1985), 50--80.

\bibitem{D91}
Xianzhe Dai, \emph{Adiabatic limits, nonmultiplicativity of signature, and
  {L}eray spectral sequence}, J. Amer. Math. Soc. \textbf{4} (1991), 265--321.

\bibitem{DF94}
Xianzhe Dai and Daniel~S. Freed, \emph{{$\eta$}-invariants and determinant
  lines}, J. Math. Phys. \textbf{35} (1994), no.~10, 5155--5194, Topology and
  physics.

\bibitem{F00}
Daniel~S. Freed, \emph{Dirac charge quantization and generalized differential
  cohomology}, Surveys in differential geometry, Surv. Differ. Geom., VII, Int.
  Press, Somerville, MA, 2000, pp.~129--194.

\bibitem{DH00}
Daniel~S. Freed and Michael Hopkins, \emph{On {R}amond-{R}amond fields and
  {$K$}-theory}, J. High Energy Phys. (2000), no.~5, Paper 44, 14.

\bibitem{FL10}
Daniel~S. Freed and John Lott, \emph{An index theorem in differential
  ${K}$-theory}, Geom. Topol. \textbf{14} (2010), 903--966.

\bibitem{GHV}
Werner Greub, Stephen Halperin, and Ray Vanstone, \emph{Connections, curvature,
  and cohomology. {V}ol. {I}: {D}e {R}ham cohomology of manifolds and vector
  bundles}, Academic Press, New York-London, 1972, {P}ure and Applied
  Mathematics, Vol. 47.

\bibitem{H14}
Man-Ho Ho, \emph{A condensed proof of the differential
  {G}rothendieck-{R}iemann-{R}och theorem}, Proc. Amer. Math. Soc. \textbf{142}
  (2014), no.~6, 1973--1982.

\bibitem{HS15}
Andreas Holmstrom and Jakob Scholbach, \emph{Arakelov motivic cohomology {I}},
  J. Algebraic Geom. \textbf{24} (2015), 719--754.

\bibitem{HS05}
M.~J. Hopkins and I.~M. Singer, \emph{Quadratic functions in geometry,
  topology,and {M}-theory}, J. Diff. Geom. \textbf{70} (2005), 329--425.

\bibitem{Kar87}
Max Karoubi, \emph{Homologie cyclique et {$K$}-th\'eorie}, Ast\'erisque (1987),
  no.~149, 147.

\bibitem{K08}
\bysame, \emph{{$K$}-theory}, Classics in Mathematics, Springer-Verlag, Berlin,
  2008, An introduction, Reprint of the 1978 edition, With a new postface by
  the author and a list of errata.

\bibitem{LM89}
H.~Blaine Lawson, Jr. and Marie-Louise Michelsohn, \emph{Spin geometry},
  Princeton Mathematical Series, vol.~38, Princeton University Press,
  Princeton, NJ, 1989.

\bibitem{L15}
Bo~Liu, \emph{Functoriality of equivariant eta forms},
  \href{http://http://arxiv.org/abs/1505.04454}{arXiv:1505.04454}.

\bibitem{L94}
John Lott, \emph{$\mathbb{R}/\mathbb{Z}$ index theory}, Comm. Anal. Geom.
  \textbf{2} (1994), 279--311.

\bibitem{MW00}
Gregory Moore and Edward Witten, \emph{Self-duality, {R}amond-{R}amond fields
  and {$K$}-theory}, J. High Energy Phys. (2000), no.~5, Paper 32, 32.

\bibitem{Q85}
Daniel Quillen, \emph{Superconnections and the {C}hern character}, Topology
  \textbf{24} (1985), no.~1, 89--95.

\bibitem{SS10}
James Simons and Dennis Sullivan, \emph{Structured vector bundles define
  differential {$K$}-theory}, Quanta of maths, Clay Math. Proc., vol.~11, Amer.
  Math. Soc., Providence, RI, 2010, pp.~579--599.

\bibitem{W85}
Edward Witten, \emph{Global gravitational anomalies}, Comm. Math. Phys.
  \textbf{100} (1985), no.~2, 197--229.

\bibitem{W98}
\bysame, \emph{D-branes and {$K$}-theory}, J. High Energy Phys. (1998), no.~12,
  Paper 19, 41 pp.\ (electronic).

\bibitem{Y75}
Zen-Ichi Yosimura, \emph{Universal coefficient sequences for cohomology
  theories of {${\rm CW}$}-spectra}, Osaka J. Math. \textbf{12} (1975), no.~2,
  305--323.

\end{thebibliography}
\end{document}